# Guido Castelnuovo
## e l'insegnamento della Matematica

*Marta Menghini*

*Abstract*

*At the beginning of the '900s Guido Castelnuovo turned his attention toward methodological, didactical, historical and applicative issues. He was an active member of the ICMI (International Commission on Mathematical Instruction) and of the Mathesis, an Italian association of mathematics teachers. He also was the author of the mathematics curricula for the new Modern Lycée and he taught the first University courses devoted to pre-service training of mathematics teachers.*
*His course on "Precision and application mathematics", inspired by Felix Klein, shows Castelnuovo's point of view on mathematics and presents a "higher point of view" on the programmes of the Modern Lyceé.*

**INTRODUZIONE**

Nei primi anni del '900 Guido Castelnuovo volse quasi completamente i suoi interessi dalla geometria algebrica verso questioni di carattere metodologico, didattico, storico e applicativo. Molto concreto fu il suo impegno nell'ambito dell'insegnamento della matematica, sia a livello istituzionale che scientifico. La sua attività si è espletata nell'*ICMI* (la commissione internazionale per l'insegnamento della matematica), nella *Mathesis* (l'associazione italiana degli insegnanti di matematica), in *ambito ministeriale* – con la stesura dei programmi del liceo moderno – e all'*Università*, con i corsi universitari dedicati ai futuri insegnanti. Tutte queste attività erano legate a una forte convinzione personale, relativa al "dovere sociale" del matematico di favorire l'apprendimento di una disciplina che è al tempo stesso "potere e gioia" (Castelnuovo, 1914), nonché a una forte concezione epistemologica, relativa al ruolo della matematica in ambito scientifico.

Il corso svolto all'Università di Roma nell'anno accademico 1913-14 su "Matematica di precisione e matematica di approssimazione", basato su una serie di lezioni già tenute da Felix Klein, aiuta a comprendere la visione di Castelnuovo sulla matematica e le sue idee sui contenuti matematici per i futuri insegnanti – collegandosi anche ai programmi del liceo moderno.

**L' ICMI (International Commission on Mathematical Instruction)**

L'ICMI è oggi un'organizzazione molto attiva, che coinvolge quasi tutti i paesi del mondo. Ha recentemente celebrato i suoi 100 anni (Menghini et al., 2008); fu infatti fondata a Roma nel 1908, durante il IV Congresso Internazionale dei Matematici, organizzato proprio da Guido Castelnuovo.
Nella quarta Sessione, relativa a *questioni filosofiche, storiche e didattiche*, venne approvata la costituzione di una commissione che svolgesse uno studio comparativo internazionale sui programmi e le metodologie d'insegnamento della matematica nelle scuole secondarie (cfr. Atti…, 1909). Il primo presidente fu Felix Klein, con George Greenhill vice-presidente e Henri Fehr segretario generale. Nel 1909 uscì il primo numero della rivista

*L'Enseignement Mathématique*, organo ufficiale dell'*ICMI* (in realtà, a quel tempo la lingua ufficiale era il francese e l'acronimo era *CIEM*).

Sempre nel 1909 furono nominati i delegati delle varie nazioni aderenti. Per l'Italia furono nominati Castelnuovo, Federigo Enriques e Giovanni Vailati. Nel 1913 Castelnuovo, insieme a Emanuel Czuber e Jacques Hadamard, verrà anche cooptato per il Comitato centrale dell'*ICMI*. (Furinghetti & Giacardi, 2008).

Le attività dell'*ICMI* riguardavano la preparazione di report di confronto su metodi e contenuti dell'insegnamento della matematica nei vari paesi e a tutti i livelli scolari, proposte per la formazione degli insegnanti, e anche proposte di riforma a livello internazionale (p.es. l'introduzione dell'analisi matematica nelle scuole superiori). Un esempio delle attività è fornito dal primo convegno ufficiale dell'ICMI, che si svolse a Milano nel settembre 1911 (Fehr, 1911). Le presentazioni erano basate sul lavoro compiuto a livello nazionale da sotto-commissioni e i temi affrontati riguardavano "il rigore e la fusione dei vari rami della matematica" e "l'insegnamento della matematica a studenti di scienze fisiche e naturali".

E' proprio Castelnuovo a illustrare la questione relativa al rigore nel confronto tra i metodi di insegnamento della geometria nei vari paesi. Nella sua sintesi Castelnuovo distingue fra:

A) *Metodo interamente logico* (Peano, Hilbert, Halsted). Tutti i postulati sono posti; si discute della loro indipendenza. Lo sviluppo è interamente logico. Non si ricorre all'intuizione; le nozioni primitive sono soggette solo alla condizione di soddisfare i postulati);

B) *Basi empiriche, svolgimento logico*. Dall'osservazione dello spazio si deducono le proposizioni primitive, su cui si fonda lo sviluppo logico. Si distinguono tre sottogruppi:

$B_A$) Tutti gli assiomi necessari sono enunciati (Sannia – D'Ovidio, Veronese, Enriques – Amaldi);

$B_B$) Una parte degli assiomi è enunciata (Euclide, Thieme);

$B_C$) Si enunciano solo gli assiomi che non hanno carattere intuitivo (Kambly, Müller);

C) *Considerazioni intuitive si alternano col metodo deduttivo* (Borel, Behrendsen-Götting). Si ricorre all'evidenza quando conviene, senza che appaia in modo preciso ciò che si ammette e ciò che si dimostra;

D) *Metodo intuitivo-sperimentale* (Perry). Si presentano teoremi come fatti che hanno carattere intuitivo o possono essere dimostrati sperimentalmente, senza che si scorga il nesso logico che li unisce.

È chiaro che la scelta $B_B$ non è stata, da parte di Euclide, di carattere didattico: la classificazione è dovuta solo al punto di vista odierno. Castelnuovo precisa che nessuna nazione adotta sistematicamente i metodi A e D; è in parte possibile seguire – sulla base dei report delle singole sotto-commissioni – lo sviluppo dei metodi di insegnamento in alcuni paesi:

ITALIA: $B_C => B_B => B_A$             FRANCIA: $B_B => B_A => C$
GERMANIA: $B_B => B_C => C$         INGHILTERRA: $B_B => B_C ( => D)$

Secondo Castelnuovo il primo passo italiano corrisponde alla riforma del 1867 (che introdusse gli elementi di Euclide come libro di testo), e fa notare che l'Italia è l'unico paese in cui si va verso un rigore maggiore. Chiariamo che non si tratta – da parte di

Castelnuovo – di una lode all'Italia. Castelnuovo ricorda anche che la maggior parte dei paesi ha unificato l'algebra con la geometria.
Dal punto di vista odierno, potremmo considerare il metodo C il metodo tipico della scuola media. In alcuni scritti Castelnuovo attribuisce un tale modo di procedere alla fisica.

**La *MATHESIS* e il *BOLLETTINO***

L'associazione tra gli insegnanti di matematica *Mathesis* fu fondata nel 1895 da Rodolfo Bettazzi. Ebbe nei primi decenni di vita un notevole ruolo culturale, grazie ai contatti stabiliti tra il mondo della scuola e il mondo universitario, e anche un ruolo politico, grazie ai rapporti diretti con il Ministero dell'Istruzione. Quello di incidere sulla politica scolastica era proprio uno degli obiettivi dell'Associazione.
I membri tenevano riunioni periodiche in cui si discutevano programmi e metodi di insegnamento, nonché le proposte che provenivano dagli ambianti ministeriali.
Inizialmente l'organo ufficiale della *Mathesis* era il *Periodico di Matematica*. Dal 1903 fu invece il *Bollettino dell'Associazione* a contenere i verbali completi delle adunanze tenute dal Comitato. Dal 1909 il bollettino prese il nome di *Bollettino della Mathesis* (che avrebbe cessato le pubblicazioni nel 1920 fondendosi con il *Periodico di Matematiche*).
Il Bollettino informava i propri lettori non solo dei contenuti delle sedute locali e nazionali dell'Associazione, in particolare dei congressi della Mathesis, ma anche dell'uscita di libri e articoli interessanti e, soprattutto, fungeva da cassa di risonanza delle iniziative estere. In particolare venivano pubblicati i report italiani elaborati per l'ICMI, e i sunti dei dibattiti internazionali già pubblicati sull'*Enseignement Mathematique*. Vi ritroviamo quindi, tra l'altro, la sintesi di Castelnuovo sull'insegnamento della geometria, di cui al paragrafo precedente (Castelnuovo, 1911), e le interessanti panoramiche sull'insegnamento della matematica nella scuola elementare (Conti, 1911), negli istituti tecnici (Scorza, 1911), nei licei (Fazzari & Scarpis, 1911).
Guido Castelnuovo presiedette la Mathesis dal 1911 al 1914, subito dopo l'anno di presidenza di Francesco Severi e prima della lunga presidenza di Federigo Enriques.

Durante la sua presidenza, Castelnuovo fu personalmente coinvolto nella politica scolastica e si dedicò, nel 1912-13, all'elaborazione dei programmi di studio e delle istruzioni per l'istituendo *ginnasio-liceo moderno*.
Quest'ultimo, istituito a titolo sperimentale come sezione del Liceo Classico nel 1911, era articolato in Ginnasio superiore (biennale) e Liceo triennale (*Legge,* 21 July 1911, n. 860). La sua istituzione avvenne sulla scia dei movimenti di riforma internazionali degli inizi del '900 che, sotto la spinta di Felix Klein, tendevano a creare una nuovo "umanesimo scientifico". Nel liceo moderno, il fine di preparare agli studi superiori, non necessariamente scientifici, si realizzava attraverso lo studio del latino, delle scienze e delle lingue moderne. La matematica interessava in quanto linguaggio adatto a descrivere i fenomeni naturali e perciò i suoi programmi (1913) presentavano significative novità. "*Il rinnovamento delle matematiche del XVII secolo è legato al rifiorire delle scienze sperimentali. In quest'ottica*", si dice "*l'insegnante dovrà far notare come i concetti fondamentali della matematica moderna, quello di funzione in particolare, siano suggeriti dalle scienze d'osservazione e, precisati*

*poi dalla matematica, abbiano a loro volta esercitato un benefico influsso sullo sviluppo di questa*" (Castelnuovo, 1912).

Il calcolo approssimato, la nozione di funzione e il calcolo infinitesimale sono introdotti per la prima volta, in modo ufficiale, proprio in questa scuola, con il suggerimento di usare un approccio sperimentale e induttivo, che completi il metodo deduttivo, e di armonizzare il corso con quello di fisica.

Le istruzioni, redatte anch'esse da Castelnuovo e pubblicate anche sul *Bollettino*, sono molto dettagliate. La novità del tipo di scuola e la metodologia pensata da Castelnuovo richiedevano evidentemente delle indicazioni molto chiare. Un esempio è fornito dal suggerimento di usare l'orario ferroviario, che a qual tempo cominciava a comparire nelle edicole, per spiegare il concetto di funzione illustrando la relazione tra spazio percorso e tempo impiegato dai vari treni. È notevole che, oltre 60 anni dopo, tale esempio sia ripreso in libri di testo di matematica innovativi (per esempio in Lombardo Radice & Mancini Proia, 1977). Torneremo sui programmi del liceo moderno nell'ultima sezione.

Dopo il periodo di presidenza, Castelnuovo continuò la sua attività nella Mathesis, a fianco di Enriques. In particolare partecipò al dibattito sulla riforma Gentile del 1923, di cui la Mathesis fu il principale interlocutore matematico (Marchi & Menghini, 2013). Il contrasto era forte, Castelnuovo rifiutò di collaborare alla stesura dei nuovi programmi e al suo posto venne incaricato Gaetano Scorza. Il nuovo ruolo della matematica, in cui si distingueva tra "*teorie*" di carattere "*algoritmico*" - che richiedono solo "*l'applicazione immediata di formule e teoremi fondamentali*" – e teorie "*che si prestano a saggiare la capacità del candidato a comprendere e far sua una rigorosa sistemazione deduttiva*", [Gentile, 1923] non si addiceva alla visione della matematica su sui si basava il liceo moderno di Castelnuovo.

**I corsi universitari rivolti ai futuri insegnanti.**

Castelnuovo si trasferì a Roma nel 1891, avendo vinto il concorso per la cattedra di Geometria analitica e proiettiva. Nel 1903, dopo la morte di Luigi Cremona, gli succedette nell'insegnamento di Geometria Superiore. Fu proprio all'interno di tale corso che, alternandoli a insegnamenti sulla geometria delle curve e superfici, Castelnuovo cominciò a impartire insegnamenti di carattere estensivo, che oggi chiameremmo di *matematiche elementari da un punto di vista superiore*. Già nell'anno accademico 1903-04 tenne un corso su *Indirizzi geometrici*, centrato sui vari tipi di trasformazioni in base alla classificazione di Felix Klein (il quale infatti inserirà questi argomenti nel secondo volume della sua serie *Elementarmathematik von einem höheren Standpunkte aus*). Castelnuovo ripeterà questo corso nel 1915-16. In anni successivi troviamo il corso sulla *Geometria non-euclidea* (1910-11 e 1919-20), il corso su *Matematica di precisione e matematica delle approssimazioni* (1913-14), il corso di *Calcolo delle probabilità* (1914-15).

A partire dall'anno accademico 1923-24 Castelnuovo iniziò gli insegnamenti di *Matematiche Complementari*, esplicitamente rivolti ai futuri insegnanti. Questi corsi nascevano sulla scia di una riflessione già iniziata da Felix Klein sui contenuti di matematici per i futuri insegnanti. L'idea era di rivedere gli argomenti che sono usualmente oggetto di insegnamento scolastico dal punto di vista delle teorie matematiche più generali e dei loro fondamenti. Klein voleva che i futuri insegnanti si rendessero conto di quanto lo studio delle matematiche superiori fosse importante proprio per la loro professione. Dallo studio universitario lo studente deve ricavare nuovi stimoli per l'insegnamento; questo è – per Klein – il principale compito dell'università.

All'inizio del corso del 1913-14 (ancora all'interno di Geometria superiore) Castelnuovo presenta le varie opinioni esistenti sulla preparazione più efficace dei futuri insegnanti, che riprenderà poi nel corso di Matematiche Complementari del 1923-24:
a) Cultura intensiva nei rami più elevati della matematica; l'attitudine didattica si formerà da se. Questa tendenza è "esagerata" quando gli uditori non abbiano l'intenzione di coltivare l'"alta scienza".
b) Cultura estensiva, allargamento della cultura nei vari indirizzi matematici e nelle scienze che con la matematica hanno la massima affinità. Questa tendenza è la più adatta per allargare le idee del futuro insegnante e per mettere nella giusta prospettiva l'argomento che si insegna. Ma non deve condurre ad un corso di carattere enciclopedico con nozioni superficiali.
c) Cultura specifica, metodologica. Questa tendenza va contenuta entro i giusti limiti, e ad ogni modo ad essa provvede la Scuola di Magistero.

Castelnuovo dichiara di intendere il corso di matematiche complementari come lo studio di questioni di matematiche superiori che hanno attinenza con le matematiche elementari. Sottolinea, come Klein, l'importanza della matematiche superiori per la formazione del futuro insegnante.
*"Anziché diffondermi sopra molti argomenti di questo tipo, tratterò con sufficiente ampiezza due argomenti: 1) Geometrie non euclidee (1° sem), 2) Problemi risolubili con riga e compasso e in particolare problema ella divisione del cerchio (2° sem).*
*Così potranno essere esaminate questioni geometriche legate col primo soggetto e questioni algebriche legate col secondo"* (Castelnuovo 1923-24)

Questa limitazione della scelta è in linea con le cautele indicate nel punto b, e spiega il perché della varietà dei corsi che saranno impartiti negli anni successivi: non un corso su tutti gli argomenti giudicati utili per i futuri insegnanti (come purtroppo abbiamo la tendenza a fare oggi), ma un corso limitato solo a un paio di argomenti, cui dedicare il necessario approfondimento. E' chiaro, in questo modo, che non sono i singoli contenuti ad essere importanti, ma il modo in cui essi vengono svolti.
I corsi di Matematiche Complementari riguardano dunque, nei vari anni:
1923-24 Geometria non euclidea. Costruzioni con riga e compasso. Poligoni regolari.
1924-25 Numeri algebrici e numeri trascendenti
1925-26 Massimi e minimi
1926-27 Lunghezze, aree, volumi dal punto di vista storico

Molto in linea con i corsi di Matematiche complementari è anche il contenuto del libro di Castelnuovo sul calcolo infinitesimale, scritto nel 1938. Non sappiamo se sarebbe diventato un corso, perché Castelnuovo fu allontanato dall'Università in seguito alle leggi razziali.

**Il corso *Matematica di precisione e matematica delle approssimazioni* (1913-14), e i programmi del liceo moderno.**

Nel 1913-14 Castelnuovo dedicò il primo semestre del Corso di Geometria superiore a "Matematica di precisione e matematica delle approssimazioni". Le lezioni erano esplicitamente ispirate a quelle tenute da Felix Klein nel 1901 e edite nel 1902 con il titolo "Anwendung der Differential- und Integralrechnung auf die Geometrie, eine Revision der Prinzipien", ripubblicate poi nel 1928 con il titolo "Präzisions- und Approximationsmathematik", come volume III della serie Matematiche elementari da un punto di vista superiore (traduzione in inglese: Klein, 2016).

Castelnuovo dedica queste lezioni ai futuri insegnanti. Dopo aver illustrato le tre opinioni sulla formazione degli insegnanti di matematica, di cui abbiamo parlato nella sezione precedente, Castelnuovo spiega perché un tale corso sia adatto alla loro formazione:

*"Alla tendenza b) può anche riallacciarsi l'opportunità di mettere bene in luce i rapporti che passano tra la matematica pura e la matematica applicata. In particolare: come si formano i concetti matematici partendo dall'osservazione del modo esterno; e come i risultati matematici possano – a loro volta – verificarsi nella realtà. Tali questioni hanno un'importanza grandissima dal punto di vista didattico. Il valore educativo della matematica sarebbe molto arricchito se, nell'insegnamento, accanto ai procedimenti logici che servono per ricavare i teoremi dai postulati, si [aggiungesse] mediante quali [procedimenti] questi si ricavano dall'osservazione, e, d'altra parte, con quale [precisione] i risultati teorici si verifichino nella realtà"*

Castelnuovo condivide con Klein non solo la convinzione dell'utilità di un tale corso per i futuri insegnanti, ma anche una visione della matematica come scienza basata sul continuo confronto tra astrazione e realtà. Nel suo articolo *Il Valore didattico della Matematica e della Fisica* (1907), che costituisce quasi un manifesto del suo pensiero, Castelnuovo dichiara che è necessario insegnare insieme matematica e fisica. Sostiene infatti l'importanza dell'osservazione, delle attività sperimentali e delle applicazioni per "*mettere in luce il valore della scienza.*" Ritiene, inoltre, che le procedure euristiche dovrebbero essere favorite per due motivi:
*"il primo, e il più importante, è che questo tipo di ragionamento è il modo migliore per giungere alla verità, non solo nelle scienze sperimentali, ma anche nella matematica stessa [...] il secondo è nel fatto che si tratta dell'unico tipo di procedimento logico che è applicabile nella vita quotidiana e in tutte le conoscenze coinvolte con essa (Castelnuovo 1907, 336)"*

In aggiunta a queste considerazioni, il corso su "Matematica di precisione e matematica di approssimazione" sembra anche motivato dalla necessità di formare i futuri insegnanti del "suo" liceo moderno, che si sta avviando in quegli anni. L'intreccio tra i programmi del corso, i programmi del liceo moderno e la visione che Castelnuovo ha della matematica è piuttosto evidente. Osserviamo che, pur essendo il corso solo scritto a mano su quaderni di appunti, contiene alcune parti molto dettagliate e diverse riflessioni. (Gario, quaderni lincei)

Il corso viene introdotto da una questione legata alla geometria (che non si trova nel testo di Klein): "Entro quali limiti i risultati della matematica teorica si trovano verificati nella realtà?". Per rispondere, Castelnuovo distingue fra tre classi di teoremi geometrici:

a) *Teoremi verificabili in modo preciso*: ad esempio la formula di Eulero per i poliedri; l'esistenza delle superfici unilatere.
b) *Teoremi verificabili in modo approssimato*: tali sono la maggior parte dei teoremi euclidei, per es.: la somma degli angoli di un triangolo è 180°; la transitività del parallelismo, ….
Castelnovo dichiara che "l'assoluta esattezza non è mai raggiungibile e non ha senso", anche usando particolari strumenti di precisione. È raggiungibile solo con gli assiomi.

In collegamento a questo punto, leggiamo nelle istruzioni ai programmi del *liceo moderno*:
*Nel ricordare agli alunni come le lunghezze e gli angoli si misurino praticamente col metro e col goniometro, l'insegnante avrà cura di avvertire che ogni misura concreta è necessariamente affetta da un errore che può essere ridotto perfezionando i mezzi di misura, ma non può mai venire soppresso.*
*Egli aggiungerà che nelle scienze applicate più evolute (geodesia, astronomia) viene prefissato un limite che l'errore non deve sorpassare e, quando tale condizione sia soddisfatta, la misura viene riguardata praticamente come esatta.*
*Le misura approssimate condurranno naturalmente l'insegnante a discorrere delle operazioni sui numeri decimali che rappresentano valori approssimati…*
*Il confronto tra le misure approssimate e le misure esatte delle grandezze fa sorgere l'idea dell'esistenza o meno di una comune misura, donde il concetto di grandezza incommensurabile. A queste si riattaccano i numeri irrazionali….*

Dunque si giunge ai numeri irrazionali partendo da misure pratiche.
E da qui riprendiamo il testo del corso universitario, che prosegue con il terzo tipo di teoremi:
c) Teoremi che non sono assolutamente verificabili: l'esistenza di segmenti incommensurabili.

Dunque, nota Castelnuovo, il ragionamento logico matematico non basta ad assicurare il risultato di un'esperienza e d'altra parte non sempre i procedimenti istintivi ed empirici sono sufficienti a giustificare i risultati matematici. È questo il motivo per cui distinguiamo fra:
*Matematica di precisione*, che riguarda tutte le proposizioni che si deducono logicamente dai postulati della geometria o dell'analisi.
*Matematica di approssimazione*: che considera i risultati ottenuti con il grado di approssimazione che l'esperienza comporta.

Dopo questa introduzione, Castelnuovo ripercorre la storia del concetto di funzione. Questa parte è più estesa della corrispondente parte nel testo di Klein e ha per scopo quello di portare gli studenti ad una definizione di funzione che comprenda tutti i casi abitualmente trattati: occorre che, dato *x*, si possa, con un numero finito di operazioni aritmetiche, ottenere *y* con un grado di approssimazione arbitrariamente grande.

Nei programmi del *Liceo Moderno* il concetto di funzione è importantissimo e vi rappresenta la vera novità:
*l'insegnante dovrà far notare come i concetti fondamentali della matematica moderna, quello di funzione in particolare, siano suggeriti dalle scienze d'osservazione e, precisati poi dalla matematica, abbiano a loro volta esercitato un benefico influsso sullo sviluppo di questa.*

La precedente affermazione contenuta nel programmi scolastici è ben esemplificata dalla *relazione tra curva empirica e curva ideale,* descritta in modo approfondito nel corso di Castelnuovo e interamente ispirata alla trattazione di Felix Klein.

Spiega infatti Castelnuovo che le successive estensioni del concetto di funzione – descritte in precedenza nel corso - hanno portato alla conseguenza che, se per curva si intende l'insieme dei punti le cui coordinate soddisfano una relazione del tipo $y = f(x)$, una tal curva può mancare di alcuni *caratteri intuitivi*. Sorge dunque la questione di esaminare quali restrizioni debbano imporsi affinché la curva ideale $y = f(x)$ presenti tali caratteri, e appaia in un certo senso come una curva empirica.

Una curva empirica (e le scienze dell'osservazione forniscono solo curve empiriche) segnata da una matita o dalla punta di uno strumento registratore, è una figura la cui larghezza non supera un certo limite, che si considera come trascurabile nelle questioni trattate. "*La curva empirica non dà dunque $y$ come funzione di $x$ nel senso della matematica di precisione, ma rappresenta solo una relazione del tipo $y = f(x) \pm \varepsilon$, dove $f$ è una funzione di $x$, ed $\varepsilon$ è un numero positivo variabile con $x$ e soggetto alla sola condizione di non superare un certo valore $\delta$. Il Klein chiama l'ente analitico $y = f(x) \pm \varepsilon$ striscia di funzione, e dice in conseguenza che una curva empirica definisce non una funzione, ma una striscia di funzione*".

Quali sono i caratteri intuitivi che si danno per scontati in una curva empirica, e quali sono quindi le condizioni affinché una curva ideale (astratta) rappresenti una curva empirica?
1) Innanzitutto la continuità.

2) Si postula poi l'esistenza di un numero finito di massimi e minimi. Situazioni come quelle rappresentate da $y = \sin \frac{1}{x}$ appartengono solo alla matematica ideale (fig.1).

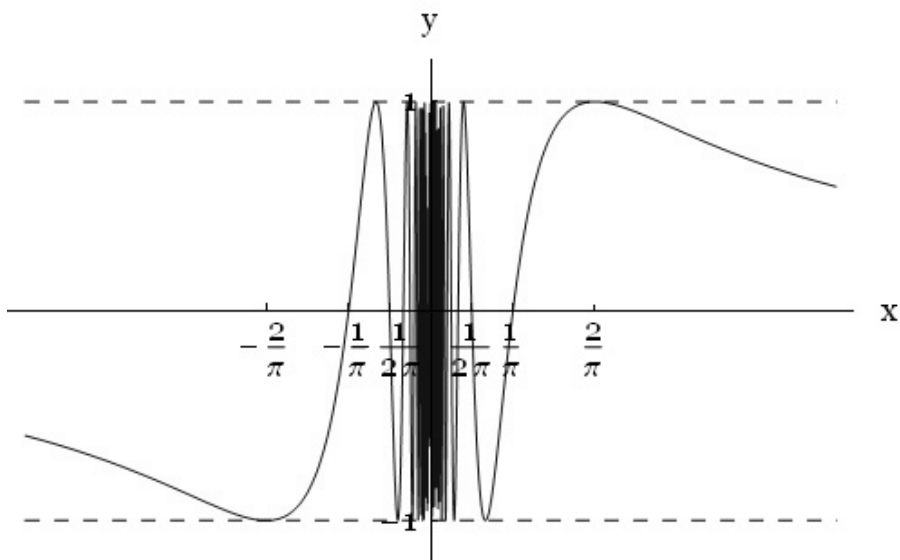

3. Derivabilità. Una curva empirica ha in ogni punto una direzione, in cui si muoverebbe un punto materiale descrivente la curva quando cessassero le forze cui è soggetto. La pendenza della curva empirica è in tal caso fornita da un rapporto $\frac{\Delta y}{\Delta x}$ (fig.2, tratta Klein, 1894), mentre nella corrispondente curva ideale chiediamo di poter calcolare la derivata in un punto.

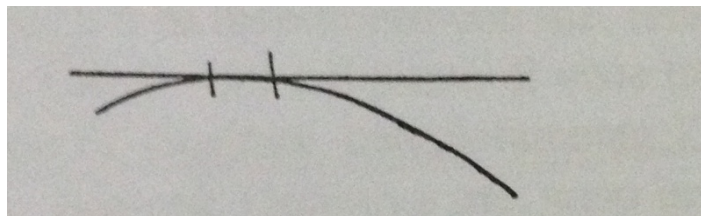

Una volta stabilite le caratteristiche di curve che *non* hanno carattere *intuitivo*, il matematico può ovviamente dedicarsi al loro studio. È una curva ideale, ad esempio, la funzione $y = \sin\frac{1}{x}$ illustrata in precedenza7. Ma, soprattutto, sia Klein che Castelnuovo dedicano molto tempo allo studio della funzione di Weierstrass (in fig. 3 le curve approssimanti disegnate nel quaderno di Castelnuovo ),

$$y = \sum_{n=0}^{n=\infty} b^n \cos(a^n \pi x)$$

che sicuramente rappresentava all'epoca una significativa novità, spiegando anche per quale motivi i matematici avessero a lungo creduto che una funzione continua fosse anche derivabile.

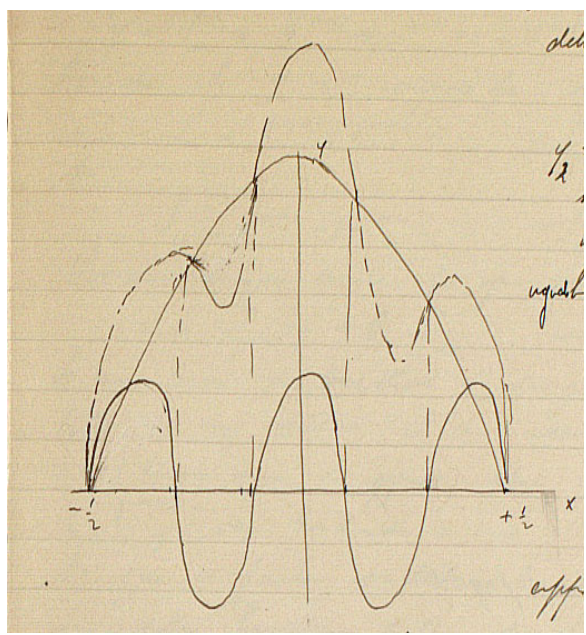

Altri argomenti, anch'essi contenuti nel testo di Klein, non sono sviluppati nel quaderno di Castelnuovo ma solo elencati. Si tratta della curva di Peano (un altro esempio di curva ideale non intuitiva), dell'interpolazione polinomiale (in particolare, Klein insiste sugli sviluppi in serie finiti e sulla loro capacità di approssimare una funzione analitica), l'analizzatore armonico. Quest'ultimo esempio è interessante. Klein descrive minuziosamente questo strumento che calcola i coefficienti di Fourier; si tratta di un apparecchio (costruito dallo Svizzero Gottlieb Coradi verso la fine dell'800) che certo non è più attuale. Tuttavia la descrizione del suo funzionamento ha un forte valore didattico e aiuta a comprendere i concetti che vi sono coinvolti.

**Conclusione**

Concludiamo semplicemente con una citazioni di Castelnuovo, che commenta bene quanto finora esposto:

*È questo il torto precipuo dello spirito dottrinario che invade la nostra scuola. Noi vi insegniamo a diffidare dell'approssimazione, che è realtà, per adottare l'idolo di una perfezione che è illusoria. […] se noi per amore della cultura soffochiamo in questi discepoli il senso pratico e lo spirito d'iniziativa, noi manchiamo al maggiore dei nostri doveri* (Castelnuovo, 1912).

## Bibliografia